\documentclass[a4paper,12pt]{article}
\usepackage{amsmath,amssymb}
\usepackage{latexsym}
\usepackage{graphicx}

\newtheorem{theorem}{Theorem}
\newtheorem{lemma}[theorem]{Lemma}

\newcommand{\qed}{~$\Box $}

\begin{document}

\title{Characterization of 1-Tough Graphs using Factors}

\author{Hongliang Lu \footnote{luhongliang@mail.xjtu.edu.cn;
Supported by the National Natural
Science Foundation of China under grant No.11471257 and
Fundamental Research Funds for the Central Universities}
\\ School of Mathematics and Statistics\\
Xi'an Jiaotong University\\
Xi'an, Shaanxi 710049, China\\
\and
Mikio Kano \footnote{mikio.kano.math@vc.ibaraki.ac.jp;
Supported by JSPS KAKENHI Grant Number 16K05248}
\\ Ibaraki University, Hitachi, Ibaraki, Japan
}

\maketitle

\begin{abstract}
For a graph $G$, let $odd(G)$ and $\omega(G)$ denote
the number of odd components and the number of
components of $G$, respectively. Then it is well-known
that $G$ has a 1-factor if and only if
$odd(G-S)\le |S|$ for all $S\subset V(G)$. Also it is clear that
$odd(G-S) \le \omega(G-S)$.
 In this paper we characterize a
1-tough graph $G$, which satisfies $\omega(G-S) \le |S|$ for all
$\emptyset \ne S \subset V(G)$, using an $H$-factor of a set-valued function
$H:V(G) \to \{ \{1\}, \{0,2\} \}$.
Moreover, we generalize this characterization to
a graph that satisfies $\omega(G-S) \le f(S)$ for all
$\emptyset \ne S \subset V(G)$,
where $f:V(G) \to \{1,3,5, \ldots\}$.
\end{abstract}

\section{Introduction}
We consider finite simple graphs, which have neither loops nor multiple edges.
Let $G$ be a graph with vertex set $V(G)$ and edge set $E(G)$.
We denote by $iso(G)$ and $odd(G)$ the number of isolated vertices
and the number of odd components of $G$, respectively.
For a set $\mathcal{S}$ of connected graphs, a spanning subgraph $F$ of $G$
is called an {\em $\mathcal{S}$-factor} if each component of $F$ is isomorphic to an element
of $\mathcal{S}$.  For an integer $n\ge 3$, let $C_n$ denote the cycle of order $n$,
and $K_2$ denote the complete graph of order 2.
Thus each component of $\{K_2,C_n: n \ge 3\}$-factor of a graph is $K_2$ or a cycle,
and a $\{K_2\}$-factor is nothing but a {\em 1-factor}.
A graph $G$ is said to be {\em factor-critical} if for every vertex $x$ of $G$,
$G-x$ has a 1-factor. We begin with the 1-factor theorem.

\begin{theorem}[The 1-factor theorem, \cite{Tutte-1} ]
A connected graph $G$  either has a 1-factor or is factor-critical  if and only if
\begin{align}
odd(G-S) \le |S| \hspace{2em} \mbox{for all}  \quad
\emptyset \ne S \subset V(G).
\label{eq-1}
\end{align}
\end{theorem}

Assume that a connected graph $G$ satisfies (\ref{eq-1}).
If $G$ has an even order,
then $G$ has a 1-factor, otherwise, $G$ is factor critical.
Moreover,  the 1-factor theorem is usually stated as follows:
a graph $G$ has a 1-factor
if and only if $odd(G-S)\le |S|$ for all $S \subset V(G)$.
By letting $S=\emptyset$ in this form,
we obtain that every component of $G$ is of even order.
However as mentioned in the above
theorem, if we use $\emptyset \ne S \subset V(G)$ instead of $S \subset V(G)$,
then the order $G$ is not necessary to be even,
and if $G$ has odd order and satisfies (\ref{eq-1}),
then $G$ is factor-critical.
This fact is shown as follows:

It is known that a graph $H$ of even order satisfies $odd(H-X)\equiv |X| \pmod{2}$
for every $X\subset V(H)$.
Assume that a connected graph $G$ has odd order
and satisfies (\ref{eq-1}), and let $x$ be any vertex of $G$.
Then $G-x$ is of even order, and for every $S\subset V(G-x)$,
it follows from (\ref{eq-1}) and the property given above that
\begin{align*}
 odd(G-x-S) & =odd(G-(S\cup \{x\})) \le |S\cup \{x\}|=|S|+1 ~~\mbox{and} \\
  odd(G-x-S) & \equiv |S| \pmod{2}.
  \end{align*}
Thus $odd(G-x-S) \le |S|$. So $G-x$ has a 1-factor by the usual 1-factor theorem,
and hence $G$ is factor-critical.
Conversely, if $G$ is factor-critical,
then for $\emptyset \ne S\subset V(G)$ and $y\in S$, we have
$odd(G-S)=odd(G-y-(S-y)) \le |S-y|\le |S|$ since $G-y$ has a 1-factor.
Hence (\ref{eq-1}) holds.

The next theorem is also well-known.
\begin{theorem}[\cite{Tutte-2}, Theorem 7.2 in \cite{AK2011}]
A connected graph $G$ of order at least $2$ has
a $\{K_2,C_n: n \ge 3\}$-factor if and only if
\begin{align}
iso(G-S) \le |S| \hspace{2em} \mbox{for all} \quad
\emptyset \ne S \subset V(G).
\label{eq-2}
\end{align}
\end{theorem}

Since $iso(G-S) \le odd(G-S)$, if a connected graph $G$ of order at least 2 satisfies (\ref{eq-1}),
then $G$ satisfies (\ref{eq-2}),
and so $G$ has a $\{K_2,C_n: n \ge 3\}$-factor. This fact is explained as follows:
Assume that $G$ satisfies
(\ref{eq-1}). If $G$ has even order, then $G$ has a 1-factor, which is a $\{K_2,C_n: n \ge 3\}$-factor.
Assume that $G$ has odd order, and let $u$ and $v$ be two adjacent vertices of $G$.
Since $G$ is factor-critical, $G-u$ has a 1-factor $M_u$
and $G-v$ has a 1-factor $M_v$.
Then $M_u\cup M_v$ is a union of two matchings of $G$, and each component of
$M_u\cup M_v$ is a $K_2$, an even cycle or a path connecting $u$ and $v$.
Hence $(M_u\cup M_v)+uv$ is
a $\{K_2,C_n: n \ge 3\}$-factor of $G$, which contains at most one odd cycle.

We denote by $\omega(G)$ the number of components of $G$.
A connected graph $G$ is said to be {\em $t$-tough}
if $|S| \ge t\omega(G-S)$ for every $S\subset V(G)$ with $\omega(G-S)>1$.
It is obvious that
\[ iso(G-S) \le odd(G-S) \le \omega(G-S) \hspace{2em}
\mbox{for all} \quad \emptyset \ne S \subset V(G). \]
In this paper, we first characterize a connected graph $G$ that satisfies $\omega(G-S)\le |S|$ for all $\emptyset \ne S \subset V(G)$.
Such a graph is called {\em 1-tough}.
Bauer,  Hakimi and Schmeichel \cite{BHF90} showed that  for any positive rational number $t$,  the $t$-tough problem, which is a problem of
checking a graph to be $t$-tough or nor,  is NP-Hard.

In this paper, we give a characterization of a 1-tough graph in terms of graph factors.
Later we generalize this characterization
by using a function $f:V(G) \to \{1,3,5, \ldots \}$.
Some results related to our theorems are found in
\cite{Ama85,Cor88,CK88,EKY2015,EJKS85,lu16,lw16+}.

\section{Characterization of 1-tough graphs}

In this section, we give a characterization of a graph $G$ that satisfies
 $\omega(G-S)\le |S|$ for all $\emptyset \ne S \subset V(G)$.
In order to state our theorem, we need some notions and definitions.
Let $\mathbf{Z}$ denote the set of integers.
For two vertices $x$ and $y$ of a graph, an edge joining $x$ to $y$ is
denoted by $xy$ or $yx$.  The degree of a vertex $v$ in a subgraph $H$ is denoted by $\deg_H(v)$.
For two vertex sets $X$ and $Y$ of $G$, not necessary to be disjoint,
we denote by
$e_G(X,Y)$ the number of edges of $G$ joining a vertex of $X$ to
a vertex of $Y$.
If $C$ is a component of $G-S$,
then we briefly write $e_G(C,S)$ for $e_G(V(C),S)$.
For a vertex set $X$ of $G$,
the subgraph  of $G$ induced by $X$ is denoted by $\langle X \rangle _G$.
For a function $h:V(G) \to \mathbf{Z}$, a subset $X\subseteq V(G)$
and a component $C$ of $G-S$ for some $S\subset V(G)$, we write
\[ h(X) := \sum_{x\in X}h(x) \quad \mbox{and} \quad
 h(C) := \sum_{x\in V(C)}h(x). \]

For any vertex $x$ of $G$, let $G^x$ denote the graph obtained from $G$
by adding a new vertex $x'$ together with a new edge $xx'$, that is, $G^x=G+xx'$.
Let $H:V(G)\rightarrow \{\{1\},\{0,2\}\}$ be a set-valued function.
So $H(v)$ is equal to $\{1\}$ or $\{0,2\}$ for each vertex $v$.
We write
\[ H^{-1}(1):=\{v\in V(G): H(v)=\{1\}\}.
\]
A spanning subgraph $F$ of $G$ is called an {\em $H$-factor} if $\deg_F(v)\in H(v)$ for
all $v\in V(G)$. This $H$-factor is also called a {\em \{1,\{0,2\}\}-factor}.
It is clear that if $G$ has an $H$-factor, then $|H^{-1}(1)|$ must be even by the Handshaking Lemma.
So  if $|H^{-1}(1)|$ is odd, then $G$ has no $H$-factor.
For a function $H:V(G)\rightarrow \{\{1\},\{0,2\}\}$ and a vertex $x$ of $G$, we define
$H^x: V(G^x) \rightarrow \{\{1\},\{0,2\}\}$ as follows.
 \begin{align}
H^x(v)=\left\{\begin{array}{ll} \{1\} &\text{if $v=x'$, }\\
H(v) &\text{otherwise}.
\end{array}\right.
\label{eq-2b}
\end{align}
A graph $G$ is said to be {\em $H$-critical} or {\em \{1,\{0,2\}\}-critical}
if $G^x$ has an $H^x$-factor for every vertex $x$ of $G$.

Let $g,f:V(G) \to \mathbf{Z}$ be functions such that $g(v) \le f(v)$ and $g(v) \equiv f(v) \pmod{2}$
for all $v\in V(G)$,
where we allow that $g(x) <0$ and $\deg_G(y)<f(y)$ for some vertices $x$ and $y$ (see Theorem 6.1 in \cite{AK2011}).
Then a spanning subgraph $F$ of $G$ is called a {\em parity $(g,f)$-factor} if
\[ g(v) \le \deg_F(v) \le f(v) \hspace{1em}\mbox{and}\hspace{1em} \deg_F(v) \equiv f(v) \pmod{2} \]
for all $v\in V(G)$.
The following theorem gives a criterion for a graph to have a parity $(g,f)$-factor.

\begin{theorem}[Lov\'asz, \cite{Lov72}, Theorem 6.1 in \cite{AK2011}]\label{lov72}
Let $G$ be a connected graph and $g,f:V(G) \to \mathbf{Z}$ such that
$g(v) \le f(v)$ and $g(v) \equiv f(v) \pmod{2}$ for all $v\in V(G)$.
Then $G$ has a parity $(g,f)$-factor if and only if for any two disjoint subsets $S,T$ of $V(G)$,
\begin{align}
\eta(S,T)=f(S)-g(T)+\sum_{x\in T}\deg_G(x) ~-e_G(S,T)-q(S,T) \geq 0,
\label{eq-17}
\end{align}
where $q(S,T)$ denotes the number of  components $C$ of $G-S-T$, called $q$-odd components, such that $f(C)+e_G(C,T)\equiv 1 \pmod 2$.
If necessary, we write $\eta(G;S,T)$ and $q(G;S,T)$ for
$\eta(S,T)$ and $q(S,T)$ to express the graph $G$.
\end{theorem}

Note that if (\ref{eq-17}) holds, then
$\eta(\emptyset,\emptyset)=-q(\emptyset,\emptyset)\ge 0$, which implies that
$|f(V(G))|\equiv 0 \pmod{2}$.
The following lemma is useful.
\begin{lemma} \label{lemma-1}
Let $G$, $g$, $f$, $S$, $T$ and  $\eta(S,T)$ be the same as Theorem~\ref{lov72}.
Then
\begin{align*}
\eta(S,T) \equiv \sum_{x\in V(G)} f(x)  \pmod{2}.
\end{align*}
\end{lemma}

\noindent {\em Proof.} Let $C_1, C_2, \ldots, C_m$ be the $q$-odd components
of $G-(S\cup T)$,
and let $D_1, D_2, \ldots,  D_r$ be the other components of  $G-(S\cup T)$. Then
$m=q(S,T)$, $f(C_i)+e_G(C_i,T)\equiv 1 \pmod{2}$ for $1\le i \le m$, and
$f(D_j)+e_G(D_j,T)\equiv 0 \pmod{2}$ for $1\le j \le r$. Hence
\begin{align*}
m  & \equiv \sum_{i=1}^m (f(C_i)+e_G(C_i,T)) + \sum_{j=1}^r (f(D_j)+e_G(D_j,T)) \\
 & \equiv \sum_{x\in V(G)-(S\cup T)} f(x) + e_G(V(G)-(S\cup T),T) \pmod{2}.
 \end{align*}
 Since $g(x)\equiv f(x) \pmod{2}$ and $-k\equiv k \pmod{2}$ for every integer $k$,
we have the following.
 \begin{align*}
 \eta(S,T) & \equiv f(S)+f(T)+ \sum_{x\in T}\deg_{G}(x) +e_G(S,T) +m \\
 & \equiv f(S) +f(T) +e_G(V(G),T) +e_G(S,T)  \\
 & \hspace{3em}  + \sum_{x\in V(G)-(S\cup T)} f(x) + e_G(V(G)-(S\cup T),T) \\
 & \equiv f(V(G)) + 2|E(\langle T \rangle _G)|   \quad (\mbox{by}~e_G(T,T) =2|E(\langle T \rangle _G)|)\\
 & \equiv \sum_{x\in V(G)} f(x)  \pmod{2}.
 \end{align*}
 Therefore the lemma holds. \qed

 \bigskip

The next  theorem is our first result, which gives a characterization
of a 1-tough graph.

\begin{theorem} \label{th-1}
Let $G$ be a connected graph. Then the following two statements hold.
\begin{description}
\item[(i)]  $G$ has an $H$-factor for every $H:V(G)\rightarrow \{\{1\},\{0,2\}\}$  with $|H^{-1}(1)|$ even if and only if
\begin{align}
\omega(G-S)\leq |S|+1\quad \mbox{for all}\quad  S\subset V(G).
\label{eq-511}
\end{align}

\item[(ii)]   $G$ is $H$-critical for every $H:V(G)\rightarrow \{\{1\},\{0,2\}\}$  with $|H^{-1}(1)|$ odd if and only if
\begin{align}
\omega(G-S)\leq |S|\quad \mbox{for all}\quad \emptyset\neq  S\subset V(G).
\label{eq-512}
\end{align}
\end{description}
\end{theorem}

\noindent {\em Proof.}
We first prove the statement (i).
Let $H:V(G)\rightarrow \{\{1\},\{0,2\}\}$ be any set-valued function
such that $|H^{-1}(1)|$ is even.
Let $M$ be a sufficiently large odd integer.
Define $f:V(G)\rightarrow \mathbf{Z}$ as
\[f(v)=\left\{
\begin{array}{ll} 1 &\text{if $H(v)=\{1\}$,}\\
2 &\text{otherwise.}
\end{array} \right.\]
Next define $g:V(G)\rightarrow \mathbf{Z}$ as
\[g(v)=\left\{
\begin{array}{ll} -M &\text{if $H(v)=\{1\}$,}\\
-M-1&\text{otherwise}.
\end{array} \right.\]

Then it is easy to see that $G$ has an $H$-factor if and only if $G$ has a parity $(g,f)$-factor.
We use Theorem \ref{lov72}. Let $S$ and $T$ be two disjoint subsets of $V(G)$.
If $T\ne \emptyset$,
then $-g(T)$ is sufficiently large, and so
\[
\eta(S,T)=f(S)-g(T)+\sum_{x\in T}\deg_G (x) -e_G(S,T) -q(S,T) \ge 0.
\]
Thus we may assume that $T=\emptyset$.
It follows that $\eta(\emptyset, \emptyset)=-q(\emptyset, \emptyset)=0$
since $f(V(G))\equiv |H^{-1}(1)|\equiv 0 \pmod{2}$ and $G$ is connected.
Hence we may assume $S\ne \emptyset$.
By $q(S,\emptyset)\le \omega(G-S)$ and (\ref{eq-511}), we have
\[
\eta(S,\emptyset)=f(S)-q(S,\emptyset )\ge |S|-\omega(G-S) \ge -1.
\]
By $f(V(G))\equiv 0 \pmod{2}$ and Lemma~\ref{lemma-1},
the above inequality implies $\eta(S,\emptyset)\ge 0$.
Therefore $G$ has the desired $H$-factor.
\medskip

We now prove the necessity.
Suppose that there exists a subset $\emptyset \ne S'\subset V(G)$ such that
\begin{align}
\omega(G-S')\geq |S'|+2.
\label{eq-7} \end{align}
Let $C_1, C_2, \ldots,C_a$ be the odd components of $G-S'$, and
$D_1, D_2, \ldots,D_b$ be the even components of $G-S'$,
where $|V(C_i)|$ is odd and $|V(D_j)|$ is even.
If $b\ge 1$, then take a vertex $w_i\in D_i$ for every $1\le i \le b$,
 and let $W\subseteq \{w_i : 1\le i \le b\}$ such that
$|W|\in \{b-1,b\}$ and $|V(G)|-|W|$ is even.
If $b=0$, then take $W\subseteq V(C_1)$ such that
$|W|\in \{0,1\}$ and $|V(G)|-|W|$ is even.

We define $H:V(G)\rightarrow \{\{1\},\{0,2\}\}$ as
\[H(v)=\left\{
\begin{array}{ll} \{0,2\} &\text{if $v\in W$, }\\
\{1\} &\text{otherwise}.
\end{array}\right.\]
Then $|H^{-1}(1)|$ is even by $H^{-1}(1)=V(G)-W$ and by the choice of $W$.
Let $M$ be a sufficiently large odd integer,
and define $f,g:V(G)\rightarrow \mathbf{Z}$ as
 \[f(v)=\left\{\begin{array}{ll} 2 &\text{if $v\in W$}\\
1 &\text{otherwise,}
\end{array}\right.\]
and
\[g(v)=\left\{\begin{array}{ll} -M-1  &\text{if $v\in W$}\\
-M &\text{otherwise}.
\end{array}\right.\]
Then it is easy to see that $G$ has an $H$-factor
if and only if $G$ has a parity $(g,f)$-factor.
We use Theorem~\ref{lov72}.
Since $f(S')= |S'|$ and $q(S', \emptyset)\ge \omega(G-S')-1$,
we obtain by (\ref{eq-7}) that
\begin{align*}
\eta(S',\emptyset) =f(S')-q(S',\emptyset) \le |S'|-\omega(G-S')+1\le -1.
\end{align*}
Therefore $G$ has no parity $(g,f)$-factor, which implies $G$ has no $H$-factor.

\medskip
We next prove the statement (ii).
Let $H:V(G)\rightarrow \{\{1\},\{0,2\}\}$ be any set-valued function
such that $|H^{-1}(1)|$ is odd.
Let $x$ be any chosen vertex of $G$, and define $H^x$ as (\ref{eq-2b}).
Then $(H^{x})^{-1}(1)=H^{-1}(1)\cup \{x'\}$ contains an even number of vertices.
We shall show that $G^x$ and $g,f$ satisfy the condition of Theorem \ref{lov72},
where $g$ and $f$ are defined as in the previous proof of the statement (i) and
$H^x(x')=\{1\}$, $f(x')=1$ and $g(x')=-M$.
Let $S$ and $T$ be two disjoint subsets of $V(G^x)=V(G)\cup \{x'\}$.
By the same argument given above, we may assume that $T=\emptyset$.
It follows that $\eta(\emptyset, \emptyset)=-q(\emptyset, \emptyset)=0$
since $f( V(G^x))\equiv |(H^{x})^{-1}(1)| \equiv 0 \pmod{2}$  and $G^x$ is connected.
Hence we may assume that $S\ne \emptyset$.
If $S$ contains $x'$, then $\omega(G^x-S)=\omega(G-(S-\{x'\})) $,
and so it follows from (\ref{eq-512}) that
\[
\eta(G^x;S,T)=f(S)-q(G^x; S,\emptyset )\ge |S|-\omega(G-(S-\{x'\})) \ge 1.
\]
Hence we may assume that $S$ does not contain $x'$.
If $S$ does not contain $x$, then
$\omega(G^x-S)=\omega(G-S) $,
and so
\[
\eta(G^x;S,T)=f(S)-q(G^x;S,\emptyset )\ge |S|-\omega(G-S) \ge 0.
\]
If $S$ contains $x$, then $\omega(G^x-S)=\omega(G-S)+1 $. Thus
\begin{align}
\eta(G^x;S,T)=f(S)-q(G^x;S,\emptyset )\ge |S|-\omega(G-S)-1 \ge -1.
\label{eq-6a}
\end{align}
On the other hand, since
\[ \sum_{v\in V(G^x)}f(v)\equiv |(H^{x})^{-1}(1)| \equiv 0 \pmod{2}, \]
it follows from Lemma~\ref{lemma-1} and (\ref{eq-6a})
that $\eta(G^x;S,T)\ge 0$.
Consequently $G^x$ has an $H^x$-factor, and
therefore $G$ is $H$-critical.
\medskip

Next we prove the necessity of (ii). Suppose that there exists a subset
$\emptyset \ne S'\subset V(G)$ such that
\begin{align}
\omega(G-S')\geq |S'|+1.
\label{eq-7b} \end{align}
Let $C_1, C_2, \ldots,C_a$ be the odd components of $G-S'$, and
$D_1, D_2, \ldots,D_b$ be the even components of $G-S'$,
where $|V(C_i)|$ is odd and $|V(D_j)|$ is even.
If $b\ge 1$, then take a vertex $w_i\in D_i$ for every $1\le i \le b$,
 and let $W\subseteq \{w_i : 1\le i \le b\}$ such that
$|W|\in \{b-1,b\}$ and $|V(G)|-|W|$ is odd.
If $b=0$, then let $W\subseteq V(C_1)$ such that
$|W|\in \{0,1\}$ and $|V(G)|-|W|$ is odd.
Moreover, choose one vertex $x$ from $S'$, and let $G^x=G+xx'$.

We define $H^x:V(G^x)\rightarrow \{\{1\},\{0,2\}\}$ as
\[H^x(v)=\left\{
\begin{array}{ll} \{0,2\} &\text{if $v\in W$, }\\
\{1\} &\text{otherwise}.
\end{array}\right.\]
Then $(H^x)^{-1}(1)=(V(G)-W)\cup \{x'\}$ and so $|(H^x)^{-1}(1)|$ is even.
Let $M$ be a sufficiently large odd integer, 
and define $f,g:V(G^x)\rightarrow \mathbf{Z}$ as
 \[f(v)=\left\{\begin{array}{ll} 2 &\text{if $v\in W$}\\
1 &\text{otherwise,}
\end{array}\right.\]
and
\[g(v)=\left\{\begin{array}{ll} -M-1  &\text{if $v\in W$}\\
-M &\text{otherwise}.
\end{array}\right.\]
Then it is easy to see that $G^x$ has an $H^x$-factor
if and only if $G^x$ has a parity $(g,f)$-factor.
We use Theorem~\ref{lov72}.
Since $f(S') = |S'|$ and $q(G^x;S',\emptyset)\ge \omega(G-S')-1+|\{x'\}|=\omega(G-S')$,
we obtain by (\ref{eq-7b}) that
\begin{align*}
\eta(G^x;S',\emptyset) =f(S')-q(G^x;S',\emptyset) \le |S'|-\omega(G-S')\le -1.
\end{align*}
Therefore $G^x$ has no parity $(g,f)$-factor, which implies $G$ is not
 $H$-critical.
Consequently, the proof of Theorem~\ref{th-1} is complete. \qed

\section{\{(1,f)-odd, even\}-factors}

In this section, we generalize Theorem~\ref{th-1} by using an odd integer
valued function $f$.
Let $G$ be a graph,  $f:V(G) \to \{1,3,5, \ldots \}$ be a function, and  
\[ 2N=\max\{f(x): x\in V(G)\}~+1\]
 be an even integer.
Define a set-valued function $H_f$ on $V(G)$ by
\begin{align} \label{eq-10}
H_f (v)= \{1,3, \ldots, f(v) \} \quad \mbox{or} \quad \{0,2, \ldots, 2N \}
\quad \mbox{for each $v\in V(G)$.}
\end{align}
Thus for a given function $f$, there are $2^{|V(G)|}$ set-valued functions $H_f$.
For a set-valued function $H_f$ on $V(G)$, define
\[ H_f^{-1}(f):= \{v\in V(G): H_f(v)=\{1,3, \ldots, f(v)\} ~\}. \]

A spanning subgraph $F$ of $G$ is called an {\em $H_f$-factor} if $\deg_F(v)\in H(v)$ for
all $v\in V(G)$. This $H_f$-factor is also called an {\em \{(1,f)-odd,even\}-factor}.
For a vertex $x$ of $G$,  we define a graph $G^x=G+xx'$. Moreover,
for a function $H_f$ on $V(G)$, define the function $H_f^x$ on $V(G^x)$ as follows.
 \[H_f^x(v)=\left\{\begin{array}{ll} \{1\} &\text{if $v=x'$,}\\
H_f(v) &\text{otherwise}.
\end{array}\right.\]
A graph is said to be {\em $H_f$-critical} or {\em \{(1,f)-odd,even\}-critical}
if $G^x$ has an $H_f^x$-factor for every vertex $x$ of $G$.

In this section, we prove the following theorem.

\begin{theorem} \label{th-2}
Let $G$ be a connected graph, and let $f:V(G) \to \{1,3,5, \ldots \}$ be
a function. Then the following two statements hold.
\begin{description}
  \item[(i)]  
$G$ has an $H_f$-factor   for every function  $H_f$ with  $|H_f^{-1}(f)|$ even if and only if
\begin{align}
\omega(G-S)\leq f(S)+1\quad \mbox{for all}\quad
S \subset V(G).
\label{eq-11}
\end{align}

 \item[(ii)] 
 $G$ is $H_f$-critical for every function $H_f$ with $|H_f^{-1}(f)|$ odd if and only if
\begin{align}
\omega(G-S)\leq f(S)\quad \mbox{for all}\quad  \emptyset\neq S\subset V(G).
\label{eq-112}
\end{align}
\end{description}

\end{theorem}

\noindent {\em Proof.} Since this theorem can be proved in a similar way
as Theorem~\ref{th-1}, we omit some details of the proof.
We first prove the sufficiency.
Assume that $G$ satisfies (\ref{eq-11}).
Let $H_f$ be any set-valued function defined by (\ref{eq-10}) such that
$|H_f^{-1}(f)|$ is even.
Let $M$ be a sufficiently large odd integer.
Define $f_1, g_1:V(G)\rightarrow \mathbf{Z}$ as
\[f_1(v)=\left\{
\begin{array}{ll} f(v) &\text{if $H_f(v)=\{1,3, \ldots, f(v)\}$,}\\
2N &\text{otherwise},
\end{array} \right.\]
and
\[g_1(v)=\left\{
\begin{array}{ll} -M &\text{if $H_f(v)=\{1,3, \ldots, f(v)\}$,}\\
-M-1 &\text{otherwise}.
\end{array} \right.\]
It is easy to see that $G$ has an $H_f$-factor
if and only if $G$ has a parity $(g_1,f_1)$-factor.
We use Theorem \ref{lov72}.
Let $S$ and $T$ be two disjoint subsets of $V(G)$.
If $T\ne \emptyset$,
then $-g_1(T)$ is sufficiently large, and so $ \eta(S,T) \ge 0$.
Thus we may assume that $T=\emptyset$.
It follows that $\eta(\emptyset, \emptyset)=-q(\emptyset, \emptyset)=0$
since $|H_f^{-1}(f)|$ is even and $G$ is connected.
Hence we may assume that $S\ne \emptyset$.
By $f_1(S) \ge f(S)$, $q(S,\emptyset)\le \omega(G-S)$
and by (\ref{eq-11}), we have
\[
\eta(S,\emptyset)=f_1(S)-q(S,\emptyset )\ge f(S)-\omega(G-S) \ge -1.
\]
Since $f_1(V(G))\equiv |H_f^{-1}(f)| \equiv 0 \pmod{2}$,
the above inequality implies $\eta(S,\emptyset) \ge 0$
by Lemma~\ref{lemma-1}.
Therefore $G$ has the desired $H_f$-factor.

\medskip
We next assume that $G$ satisfies (\ref{eq-112}).
Let $x$ be any chosen vertex of $G$.
We shall show that $G^x$ and $f_1, g_1$ satisfies the condition of
Theorem \ref{lov72},
where $f_1(x')=1$ and $g_1(x')=-M$.
Let $S$ and $T$ be two disjoint subsets of $V(G^x)=V(G)\cup \{x'\}$.
By the same argument given above, we may assume $T=\emptyset$.
It follows that $\eta(G^x;\emptyset, \emptyset)=-q(G^x;\emptyset, \emptyset)=0$
since
$\{v\in V(G^x): f_1(v)\equiv 1 \pmod{2}\}=\{x'\}\cup H_f^{-1}(f)$
contains an even number of vertices and $G^x$ is connected.
Hence we may assume that $S\ne \emptyset$.
If $S$ contains $x'$, then $\omega(G^x-S)\le\omega(G-(S-x')) $,
and so $ \eta(G^x;S,T) \ge f(S)-\omega(G-(S-x')) \ge 1$.
Thus we may assume that $S$ does not contain $x'$.
If $S$ does not contain $x$, then
$\omega(G^x-S)=\omega(G-S) $,
and so $\eta(G^x;S,T) \ge f(S)-\omega(G-S) \ge 0$.
If $S$ contains $x$, then $\omega(G^x-S)=\omega(G-S)+1 $, and thus
$\eta(G^x;S,T) \ge f(S)-\omega(G-S)-1\ge -1$,
which implies $\eta(G^x;S,T)\ge 0$ by Lemma~\ref{lemma-1}
and $f_1(V(G^x)) \equiv |H_f^{-1}(f)\cup \{x'\}| \equiv 0 \pmod{2}$.
Therefore $G^x$ has a $H_f^x$-factor.
Consequently $G$ is $H_f$-critical.

\bigskip

We now prove the necessity.
First consider (i).  Assume that there exists a subset
$\emptyset \ne S'\subset V(G)$ such that
\begin{align}
\omega(G-S')\geq f(S')+2.
\label{eq-7c} \end{align}
Let $C_1, C_2, \ldots,C_a$ be the odd components of $G-S'$, and
$D_1, D_2, \ldots,D_b$ be the even components of $G-S'$.
If $b\ge 1$, then take a vertex $w_i\in D_i$ for every $1\le i \le b$,
 and let $W\subseteq \{w_i : 1\le i \le b\}$ such that
$|W|\in \{b-1,b\}$ and $|V(G)|-|W|$ is even.
If $b=0$, then take $W\subseteq V(C_1)$ such that
$|W|\in \{0,1\}$ and $|V(G)|-|W|$ is even.

We define $H_f:V(G)\rightarrow \{\{1,3, \ldots, f(v)\},\{0,2, \ldots, 2N\}\}$ as
\[H_f(v)=\left\{
\begin{array}{ll} \{0,2,\ldots, 2N\} &\text{if $v\in W$, }\\
\{1,3, \ldots, f(v)\} &\text{otherwise}.
\end{array}\right.\]
Then $|H_f^{-1}(f)|$ is even by $H_f^{-1}(f)=V(G)-W$
and by the choice of $W$.
Let $M$ be a sufficiently large odd integer, 
and define $f_2,g_2:V(G)\rightarrow \mathbf{Z}$ as
 \[f_2(v)=\left\{\begin{array}{ll} 2N &\text{if $v\in W$}\\
f(v) &\text{otherwise,}
\end{array}\right.\]
and
\[g_2(v)=\left\{\begin{array}{ll} -M-1  &\text{if $v\in W$}\\
-M &\text{otherwise}.
\end{array}\right.\]
Then $G$ has an $H_f$-factor
if and only if $G$ has a parity $(g_2,f_2)$-factor.
We use Theorem~\ref{lov72}.
Since $f_2(S')= f(S')$ and $q(S', \emptyset)\ge \omega(G-S')-1$,
it follows from (\ref{eq-7c}) that
\begin{align*}
\eta(S',\emptyset) =f_2(S')-q(S',\emptyset) \le f(S')-\omega(G-S')+1\le -1.
\end{align*}
Therefore $G$ has no $H_f$-factor.

\medskip

Next consider (ii). Suppose that there exists a subset
$\emptyset \ne S'\subset V(G)$ such that
\begin{align}
\omega(G-S')\geq f(S')+1.
\label{eq-7d} \end{align}
Let $C_1, C_2, \ldots,C_a$ be the odd components of $G-S'$, and
$D_1, D_2, \ldots,D_b$ be the even components of $G-S'$.
If $b\ge 1$, then take a vertex $w_i\in D_i$ for every $1\le i \le b$,
 and let $W\subseteq \{w_i : 1\le i \le b\}$ such that
$|W|\in \{b-1,b\}$ and $|V(G)|-|W|$ is odd.
If $b=0$, then let $W\subseteq V(C_1)$ such that
$|W|\in \{0,1\}$ and $|V(G)|-|W|$ is odd.
Define a set-valued function $H_f$ on $V(G)$ as
\[H_f(v)=\left\{
\begin{array}{ll} \{0,2,\ldots, 2N\} &\text{if $v\in W$, }\\
\{1,3, \ldots, f(v)\} &\text{otherwise.}
\end{array}\right.\]
Then $|(H_f)^{-1}(f)|=|V(G)-W|$ is odd.

Choose one vertex $x$ from $S'$, and let $G^x=G+xx'$.
Then define a function$H_f^x$ on $V(G^x)$ as given before.
Let $M$ be a sufficiently large odd integer, 
and define $f_2,g_2:V(G^x)\rightarrow \mathbf{Z}$ as
 \[f_2(v)=\left\{\begin{array}{ll} 2N &\text{if $v\in W$,}\\
f(v) &\text{if $v\in V(G)-W$,} \\
1 &\text{if $v=x'$.}
\end{array}\right.\]
and
\[g_2(v)=\left\{\begin{array}{ll} -M-1  &\text{if $v\in W$}\\
-M &\text{if otherwise.}
\end{array}\right.\]
Then it is easy to see that $G^x$ has an $H_f^x$-factor
if and only if $G^x$ has a parity $(g_2,f_2)$-factor.
We use Theorem~\ref{lov72}.
Since $f_2(S') = f(S')$ and $q(G^x;S',\emptyset)\ge \omega(G-S')-1+|\{x'\}|=\omega(G-S')$,
we obtain by (\ref{eq-7d}) that
\begin{align*}
\eta(G^x;S',\emptyset) =f(S')-q(G^x;S',\emptyset) \le f(S')-\omega(G-S')\le -1.
\end{align*}
Therefore $G^x$ has no parity $(g_2,f_2)$-factor, which implies $G$ is not
 $H_f$-critical.
Consequently, the proof of Theorem~\ref{th-2} is complete. \qed

\medskip \noindent
{\bf Acknowledgment} The authors would like to thank Dr. Kenta Ozeki for his valuable suggestions and comments.

\end{document}